\documentclass[graybox]{svmult}

\usepackage{type1cm}        
\usepackage{makeidx}         
\usepackage{graphicx}        
\usepackage{multicol}        
\usepackage[bottom]{footmisc}
\usepackage[numbers]{natbib}
\usepackage{orcidlink}

\usepackage{newtxtext}       %
\usepackage[varvw]{newtxmath}       

\usepackage{amsmath,graphicx,textcomp,hyperref,tabularx,xcolor,theorem}

\hypersetup{
	colorlinks=true, 
	linkcolor=blue, 
	urlcolor=red, 
	citecolor=magenta,
	linktoc=all 
}

\usepackage{algorithmicx}
\usepackage{algorithm}
\usepackage{algpseudocode}

\usepackage{xpatch}
\makeatletter
\xpatchcmd{\algorithmic}
{\ALG@tlm\z@}{\leftmargin\z@\ALG@tlm\z@}
{}{}
\makeatother

\newcommand{\assign}{:=}
\newcommand{\cdummy}{\cdot}
\newcommand{\mathLaplace}{\Delta}

\newcommand{\tmmathbf}[1]{\ensuremath{\boldsymbol{#1}}}
\newcommand{\tmop}[1]{\ensuremath{\operatorname{#1}}}
\newcommand{\tmrsub}[1]{\ensuremath{_{\textrm{#1}}}}

\newcommand{\tmtextit}[1]{\text{{\itshape{#1}}}}

\newcommand{\tmverbatim}[1]{\text{{\ttfamily{#1}}}}

\newenvironment{tmparmod}[3]{\begin{list}{}{\setlength{\topsep}{0pt}\setlength{\leftmargin}{#1}\setlength{\rightmargin}{#2}\setlength{\parindent}{#3}\setlength{\listparindent}{\parindent}\setlength{\itemindent}{\parindent}\setlength{\parsep}{\parskip}} \item[]}{\end{list}}

\newcounter{tmcounter}

\usepackage[T3,T1]{fontenc}
\DeclareSymbolFont{tipa}{T3}{cmr}{m}{n}
\DeclareMathAccent{\invbreve}{\mathalpha}{tipa}{16}

\newcommand{\avg}[1]{\overline{#1}}
\newcommand{\myvector}[1]{\text{{\hspace{0.1em}}#1{\hspace{0.15em}}}}

\newcommand{\vu}{\myvector{u}}

\newcommand{\vD}{\myvector{D}}

\newcommand{\vf}{\myvector{f}}
\newcommand{\vg}{{\myvector{g}}}
\newcommand{\vs}{\myvector{s}}
\newcommand{\vF}{\myvector{F}}

\newcommand{\vS}{\myvector{S}}

\newcommand{\bp}{\tmmathbf{p}}

\newcommand{\bzero}{\tmmathbf{0}}
\newcommand{\bS}{\tmmathbf{S}}

\newcommand{\bx}{\tmmathbf{x}}

\newcommand{\emh}{e - \frac{1}{2}}

\newcommand{\eph}{e + \frac{1}{2}}

\newcommand{\Nmh}{N - \frac{1}{2}}

\newcommand{\pd}[2]{\frac{\partial #1}{\partial #2}}
\newcommand{\od}[2]{\frac{\ud #1}{\ud #2}}

\newcommand{\Uad}{\mathcal{U}_{\textrm{\tmop{ad}}}}

\newcommand{\uu}{\tmmathbf{u}}

\newcommand{\bw}{\tmmathbf{u}}

\newcommand{\atu}{\invbreve{\tmmathbf{u}}}
\newcommand{\utilow}{\invbreve{\tmmathbf{u}}^{\tmop{low}, n + 1}}
\newcommand{\au}{\avg{\uu}}

\newcommand{\pf}{\tmmathbf{f}}

\newcommand{\ff}{\pf}

\newcommand{\F}{\tmmathbf{F}}

\newcommand{\bss}{\tmmathbf{s}}

\newcommand{\uep}{\tmmathbf{u}_{e, p}}
\newcommand{\ueppone}{\tmmathbf{u}_{e, p + 1}}

\newcommand{\uez}{\tmmathbf{u}_{e, 0}}

\newcommand{\uepoz}{\tmmathbf{u}_{e + 1, 0}}
\newcommand{\ueN}{\tmmathbf{u}_{e, N}}

\newcommand{\pph}{p + \frac{1}{2}}
\newcommand{\pmh}{p - \frac{1}{2}}

\newcommand{\paragraphtoc}[1]{}
\newcommand{\ad}{P}

\providecommand{\subindex}{=}

\providecommand{\citep}{}
\providecommand{\citet}{}

\newcommand{\re}{\mathbb{R}}

\newcommand{\trace}{\tmop{trace}}
\newcommand{\poly}{\mathbb{P}}

\newcommand{\ud}{d}

\newcommand{\half}{\frac{1}{2}}

\newcommand{\cE}{E}

\newcommand{\Eonetwo}{E_{12}}
\newcommand{\Etwotwo}{E_{22}}
\newcommand{\vone}{v_1}
\newcommand{\vtwo}{v_2}
\newcommand{\ccE}{\tmmathbf{E}}

\newcommand{\cR}{p}
\newcommand{\ccP}{\tmmathbf{p}}
\newcommand{\ccR}{\ensuremath{\tmmathbf{p}}}
\newcommand{\Poneone}{p_{11}}
\newcommand{\Ponetwo}{p_{12}}
\newcommand{\Ptwotwo}{p_{22}}
\newcommand{\bv}{\tmmathbf{v}}

\usepackage{url}

\makeindex             

\begin{document}
	
\title*{Generalized framework for admissibility preserving Lax-Wendroff Flux
	Reconstruction for hyperbolic conservation laws with source terms}
\titlerunning{Admissibility preserving LWFR with sources}
\author{Arpit Babbar\orcidID{0000-0002-9453-370X} and\\ Praveen Chandrashekar\orcidID{0000-0003-1903-4107}}
\authorrunning{Babbar, Chandrashekar}
\institute{Arpit Babbar \at TIFR-CAM, Bangalore - 560065, \email{arpit@tifrbng.res.in}
\and Praveen Chandrashekar \at TIFR-CAM, Bangalore - 560065, \email{praveen@math.tifrbng.res.in}}
\maketitle
\abstract*{Lax-Wendroff Flux Reconstruction (LWFR) is a single-stage, high order,
	quadrature free method for solving hyperbolic conservation laws. We perform a cell average decomposition of the LWFR scheme that is similar to the one used in the admissibility preserving framework of Zhang and Shu (2010). By performing a flux limiting of the time averaged numerical flux, the decomposition is used to obtain an admissibility preserving LWFR scheme. The admissibility preservation framework is further extended to a newly proposed extension of LWFR scheme for conservation laws with source terms. This is the first extension of the high order LW scheme that can handle source terms. The admissibility and accuracy are verified by numerical experiments on the Ten Moment equations of Livermore et al.}
\section{Introduction}
Lax-Wendroff (LW) methods are a class of single stage explicit high
order methods for solving time dependent PDEs, and are an alternative to the
PDE solvers where multistage explicit Runge-Kutta (RK) methods are used for
temporal discretization. The single step nature of LW methods decreases
inter-element communication and makes them suitable for modern memory
bandwidth limited hardware. Arbitrary high order of accuracy is achieved in
Lax-Wendroff schemes by performing a high order Taylor's expansion of the
solution and using the PDE to replace the temporal derivatives with a spatial derivative of the time-averaged flux
approximation; the procedure is a generalization of~{\cite{Lax1960}}. Some
of the pioneering works are
{\cite{Qiu2003,Qiu2005b}}. Another class of single stage methods is ADER, where a high order element local solution is obtained by solving a space-time implicit equation~{\cite{Titarev2002,Dumbser2008}}.

Flux Reconstruction (FR) is a collocation based, Spectral Element
Method introduced by Huynh~{\cite{Huynh2007}} and uses the same polynomial basis
approximation as the Discontinuous Galerkin (DG) method of Cockburn and
Shu~{\cite{cockburn2000}}. The FR scheme is quadrature free and can
generalize various numerical schemes including variants of DG by making
different choices of correction functions and solution points~{\cite{Huynh2007,Trojak2021}}. The
accuracy and stability of FR has been studied
in~{\cite{Huynh2007,Vincent2011a,Trojak2021,Cicchino2022a}}.

LW schemes with FR spatial approximation were first introduced
in~{\cite{Lou2020}}, the role of FR in these schemes is to correct the time averaged
flux while coupling elements across interfaces. In~{\cite{babbar2022}}, the
present authors proposed a Lax-Wendroff Flux Reconstruction (LWFR) scheme
using the Jacobian-free approximate Lax-Wendroff procedure
of~{\cite{Zorio2017,Burger2017}}. The numerical flux was carefully
constructed in~{\cite{babbar2022}} to obtain enhanced accuracy and
stability. In~{\cite{babbar2023admissibility}}, a subcell based shock
capturing blending scheme was introduced for LWFR based on~{\cite{hennemann2021}}.~{\cite{babbar2023admissibility}} exploited the subcell structure to construct a
 \tmtextit{blended numerical flux} between the time averaged flux and the low order flux, which was used to obtain provably admissibility preserving property of the Lax-Wendroff scheme.

In this work, it is shown that the construction of \tmtextit{blended numerical flux} to enforce admissibility can also be done when there is no subcell based limiting scheme. The initial argument is similar to performing a decomposition of the cell average into \textit{fictitious finite volume updates} as in~\cite{Zhang2010b,zhang2010c}. The difference from~\cite{Zhang2010b} arises in the case of LW schemes as some of the fictitious finite volume updates involve the LW high order fluxes. Then, it is seen that, if the LW numerical flux is limited to ensure that the update with its fictitious finite volume fluxes is admissible, the scheme will be admissibility preserving in means. The limiting procedure
of~\cite{babbar2023admissibility} is then used to enforce admissibility.

To numerically validate our claim, we test LWFR on the
Ten Moment equations, which are derived by Levermore et al.~{\cite{Levermore_1996}} by taking a Gaussian closure of the kinetic model. These equations are relevant in many
applications, especially related to plasma
flows~(see~{\cite{Berthon_TMP_2006,Berthon2015}} and further
references in~{\cite{meena2017}}), in cases where the \tmtextit{local
	thermodynamic equilibrium} used to close the Euler equations of compressible
flows is not valid, and anisotropic nature of the pressure needs to be
accounted for. The Ten Moment equations model is a hyperbolic conservation
law with source terms. Thus, in addition to showing that our positivity
preserving framework preserves admissibility in the presence of shocks and
rarefactions, we also introduce the first LWFR scheme in the presence of
source terms. The approach involves adding time averages of the sources and thus we also propose a source term limiting procedure so that admissibility is still maintained.

The rest of the paper is organized as follows. Section~\ref{sec:lwfr}
introduces the LWFR scheme with source terms and reviews notions of admissibility preservation.
Section~\ref{sec:adm.pres} describes the additional limiting required in LW scheme for admissibility preservation, i.e., for the time averaged flux~(Section~\ref{sec:flux.correction}) and time averaged sources~(Section~\ref{sec:source.limiter}).
Section~\ref{sec:numerical.results} shows the numerical results for the Ten
Moment equations model and conclusions of the work are drawn in Section~\ref{sec:conclusion}.
\section{Lax-Wendroff Flux Reconstruction}\label{sec:lwfr}
Consider a conservation law of the form
\begin{equation}
	\uu_t + \ff_x = \bss \label{eq:con.law}
\end{equation}
where $\uu \in \re^p$ is the vector of conserved quantities, $\ff = \pf (\bw)$
is the corresponding flux, $\bss = \bss \left( \uu, t, x \right)$ is the
source term, together with some initial and boundary conditions. The solution
that is physically correct is assumed to belong to an admissibility set,
denoted by $\Uad$. For example in case of compressible flows, the density and
pressure (or internal energy) must remain positive. In case of the Ten Moment
equations~\eqref{eq:tmp}, the density must remain positive and the pressure
tensor $\ccR$ must be positive definite. In both these models, and most of the models that are of
interest, the admissibility set is a convex subset of $\re^p$, and can be
written as
\begin{equation}
	\label{eq:uad.form} \Uad = \{ \uu \in \re^p : \ad_k (\bw) > 0, 1 \le k \le
	K\}
\end{equation}
For Euler's equations, $K = 2$ and $\ad_1, \ad_2$ are density, pressure
functions respectively; if the density is positive then pressure is a concave
function of the conserved variables. For the Ten Moment
equations~\eqref{eq:tmp}, $K = 3$ and $P_1, P_2, P_3$ are density, $\trace
\left( \bp \right)$, $\det \left( \ccR \right)$. Although density and trace
functions are concave functions of the conserved variables, $\det \left( \ccR
\right)$ is not so.

For the numerical solution, we will divide the computational domain $\Omega$
into disjoint elements $\Omega_e$, with $\Omega_e = [x_{\emh}, x_{\eph}]$ so that $\Delta x_e =
x_{\eph} - x_{\emh}$. Then, we map each element $\Omega_e$ to the reference element $[0, 1]$ by $x \mapsto  \frac{x - x_{e-1/2}}{\Delta x_e} =: \xi$. Inside each element, we approximate the solution by degree $N \ge 0$
polynomials belonging to the set $\poly_N$. For this, choose $N + 1$ distinct
nodes $\{ \xi_i \}_{i \subindex 0}^N$ in $[0,1]$ which will be taken to be Gauss-Legendre (GL) points in this work, and will also be referred to as \tmtextit{solution points}. The solution inside an element $\Omega_e$ is given by $\bw_h  (\xi, t) = \sum_{p = 0}^N \uep (t) \ell_p(\xi)$ where $\{\ell_p\}$ are Lagrange polynomials of degree $N$ defined to satisfy $\ell_p (\xi_q) = \delta_{pq}$. Note that the coefficients $\uep$ which are the basic unknowns or
\tmtextit{degrees of freedom} (dof), are the solution values at the solution
points.

The Lax-Wendroff scheme combines the spatial and temporal discretization into a single step. The starting point is a Taylor's expansion in time following the Cauchy-Kowalewski procedure where the PDE is used to rewrite some of the temporal derivatives in the Taylor expansion as spatial derivatives. Using Taylor's expansion in time around $t = t_n$, we can write the solution at the next time
level as
\[ 
\uu^{n + 1} = \uu^n + \sum_{m = 1}^{N + 1} \frac{\Delta t^m}{m!} \partial_t^m  \uu^n + O (\Delta t^{N + 2})
\]
Since the spatial error is expected to be of $O (\Delta x^{N + 1})$, we retain
terms up to $O (\Delta t^{N + 1})$ in the Taylor expansion, so that the overall formal accuracy is of order $N + 1$ in both space and time. Using the PDE, $\partial_t  \uu = - \partial_x  \pf + \bss$, we re-write time derivatives of
the solution in terms of spatial derivatives of the flux and source terms
\[
\partial_t^m  \uu = - (\partial_t^{m - 1}  \pf)_x + \partial_t^{m - 1} \bss, \qquad m = 1, 2, \ldots 
\]
so that
\begin{equation}
\begin{split}
\uu^{n + 1} & = \uu^n - \sum_{m = 1}^{N + 1} \frac{\Delta t^m}{m!} 
(\partial_t^{m - 1}  \pf)_x + \sum_{m = 1}^{N + 1} \frac{\Delta t^m}{m!}
\partial_t^{m - 1}  \bss + O (\Delta t^{N + 2})\\
& = \uu^n - \Delta t \pd{\F}{x} (\uu^n) + \mathLaplace t \bS \left(
\uu^n, t^n \right) + O (\Delta t^{N + 2})
\end{split}
\label{eq:lwtay}
\end{equation}
where
\begin{align}
\label{eq:tavgflux} \F &= \sum_{m = 0}^N \frac{\Delta t^m}{(m + 1) !}
\partial_t^m  \pf = \pf + \frac{\Delta t}{2} \partial_t  \pf
+ \ldots + \frac{\Delta t^N}{(N + 1) !} \partial_t^N  \pf \\
\label{eq:tavgS} \bS &= \sum_{m = 0}^N \frac{\Delta t^m}{(m + 1) !}
\partial_t^m  \bss = \bss + \frac{\Delta t}{2} \partial_t  \bss
+ \ldots + \frac{\Delta t^N}{(N + 1) !} \partial_t^N  \bss
\end{align}
Note that $\F ({\uu^n}), \bS ( \uu^n, t^n)$ are approximations
to the time average flux and source term in the interval $[t_n, t_{n + 1}]$
since they can be written as
\begin{align}
\F ({\uu^n}) &= \frac{1}{\Delta t}  \int_{t_n}^{t_{n + 1}} \left[
\pf (\uu^n) + \ldots + \frac{(t -
t_n)^N}{N!} \partial_t^N  \pf (\uu^n) \right]  \ud t \label{eq:tvgproperty} \\
\bS ({\uu^n}, t^n) &= \frac{1}{\Delta t}  \int_{t_n}^{t_{n + 1}} \left[
\bss (\uu^n, t^n) + \ldots +
\frac{(t - t_n)^N}{N!} \partial_t^N  \bss (\uu^n, t^n) \right]  \ud t
\label{eq:tvgproperty.S}
\end{align}
where the quantity inside the square brackets is the truncated Taylor
expansion of the flux $\pf$ or source $\bss$ in time.
Following equation~\eqref{eq:lwtay} we need to specify the construction of the time averaged flux~\eqref{eq:tavgflux} and the time averaged source terms~\eqref{eq:tavgS}. The first step of approximating~\eqref{eq:lwtay} is the predictor step where a locally degree $N$ approximation $\F^\delta$ of the time averaged flux is computed by the approximate Lax-Wendroff procedure of Zorio~\cite{Zorio2017}, also discussed in Section 4.4 of~\cite{babbar2022}. Then, as in the standard RKFR scheme, we perform the Flux Reconstruction procedure on $\F^\delta$ to construct a locally degree $N+1$ and globally continuous flux approximation $\F_h (\xi)$. The time average source $\bS$ will also be approximated locally as a degree $N$ polynomial using the approximate Lax-Wendroff procedure and denoted  with a single notation $\bS^\delta (\xi)$ since it needs no correction. The scheme for local approximation is discussed in Section~\ref{sec:approximate.lw}. After computing $\F_h, \bS^\delta$, truncating equation~\eqref{eq:lwtay}, the solution at the nodes is updated by a collocation scheme as follows
\begin{equation}
	\uep^{n + 1} = \uep^n - \frac{\Delta t}{\Delta x_e}  \od{\F_h}{\xi} (\xi_p)
	+ \mathLaplace t \bS^\delta (\xi_p), \qquad 0 \le p \le N \label{eq:uplwfr}
\end{equation}
This is the single step Lax-Wendroff update scheme for any order of accuracy.
\subsection{Approximate Lax-Wendroff procedure for degree $N = 2$} \label{sec:approximate.lw}
The approximations of temporal derivatives of $\bss$ are made in a similar
fashion as those of $\pf$ in~\cite{Zorio2017,babbar2022}. For example, to obtain second order accuracy, $\partial_t  \bss$ can be approximated as
\begin{align*}
	\partial_t  \bss \left( \uu, \bx, t \right)  & \approx \frac{\bss \left( \uu^n + \mathLaplace t \uu^n_t, \bx, t^{n
			+ 1} \right) - \bss \left( \uu^n - \mathLaplace t \uu_t^n, \bx, t^{n - 1}
		\right)}{2 \mathLaplace t}
\end{align*}
where $\uu_t = - \partial_x  \pf + \bss \left( \uu, \bx, t \right)$. Denoting $\vg^{(k)}$ as an approximation for $\Delta t^k \partial_t^k g$, we explain the local flux and source term approximation procedure for degree $N = 2$
\begin{align*}
	\vF = \vf + \frac{1}{2} \vf^{(1)} + \frac{1}{6} \vf^{(2)}, \qquad \vS = \vs + \frac{1}{2}  \vs^{(1)} + \frac{1}{6}  \vs^{(2)}
\end{align*}
where
\begin{eqnarray*}
	\vu^{(1)} &=& - \frac{\Delta t}{\Delta x_e} \vD \vf + \Delta t \vs\\
	\vf^{(1)}, \vs^{(1)} &=& \frac{1}{2} \left[ \pf\left(\vu + \vu^{(1)}\right) - \pf\left(\vu - \vu^{(1)}\right) \right], \frac{1}{2} \left[ \bss\left(\vu + \vu^{(1)}\right) - \bss\left(\vu - \vu^{(1)}\right) \right] \\
	\vu^{(2)} &=& - \frac{\Delta t}{\Delta x_e} \vD \vf^{(1)}  + \Delta t \vs^{(1)} \\
	\vf^{(2)} &=& \pf\left(\vu + \vu^{(1)} + \half  \vu^{(2)}\right) - 2 \pf(\vu) + \pf\left(\vu -  \vu^{(1)} + \half  \vu^{(2)} \right) \\
	\vs^{(2)} &=& \bss\left(\vu + \vu^{(1)} + \half  \vu^{(2)}\right) - 2 \bss(\vu) + \bss\left(\vu -  \vu^{(1)} + \half  \vu^{(2)} \right)
\end{eqnarray*}
For complete details on the local approximation of the flux $\F$ for all
degrees and then its FR correction using the time numerical flux $\F_{\eph}$,
the reader is referred to~{\cite{babbar2022}}.
\subsection{Admissibility preservation}
The admissibility preserving property of the conservation law, also known as convex set preservation
property since $\Uad$ is convex, can be written as
\begin{equation}
	\label{eq:conv.pres.con.law} \bw (\cdummy, t_0) \in \Uad \qquad
	\Longrightarrow \qquad \bw (\cdummy, t) \in \Uad, \qquad t > t_0
\end{equation}
Thus, an admissibility preserving FR scheme is one which preserves admissibility at all solution points (Definition 1 of~\cite{babbar2023admissibility}). In this work, as in~\cite{babbar2023admissibility}, we study the admissibility preservation in means property of the LWFR scheme (Definition 2 of~\cite{babbar2023admissibility}) which is said to hold for schemes that satisfy
\begin{equation} \label{defn:mean.pres}
	\uep^n \in \Uad \quad \forall e, p \qquad \Longrightarrow \qquad \au_e^{n + 1} \in \Uad \quad \forall e 
\end{equation}
where $\avg{\uu}_e := \sum_{p=0}^N w_p \uu_{e,p}$ denotes the cell average of $\uu_e$. Once the scheme is admissibility preserving in means, the scaling limiter of~{\cite{Zhang2010b}} can be used to
obtain an admissibility preserving scheme. The following property of the LWFR scheme will be crucial in the obtaining admissibility preservation in means
\begin{equation}
\label{eq:fravgup} \overline{\bw}_e^{n + 1} = \au_e^n - \frac{\Delta t}{\Delta x_e}  (\F_{\eph} - \F_{\emh}) + \Delta t \avg{\bS}_e
\end{equation}
where $\avg{\bS}_e:= \sum_{p=0}^N w_p \bS_e^\delta(\xi_p)$ is the cell average of the source term.
\section{Limiting time averages} \label{sec:adm.pres}
\subsection{Limiting time average flux} \label{sec:flux.correction}
In this section, we study admissibility preservation in means property~\eqref{defn:mean.pres} for the LWFR update~\eqref{eq:uplwfr} in the case where source term $\bss$ in~\eqref{eq:con.law} is zero. Similar to the work of Zhang-Shu~{\cite{Zhang2010b}}, we define \textit{fictitious finite volume updates}
\begin{equation}
	\begin{split}
		\atu_{e, 0}^{n + 1} & = \uez^n - \frac{\Delta t}{w_0 \Delta x_e} 
		[\pf_{\half}^e - \F_{\emh}^{\tmop{LW}}]\\
		\atu_{e, p}^{n + 1} & = \uep^n - \frac{\Delta t}{w_p \Delta x_e} 
		[\pf_{\pph}^e - \pf_{\pmh}^e], \qquad 1 \le p \le N - 1\\
		\atu_{e, N}^{n + 1} & = \ueN^n - \frac{\Delta t}{w_N \Delta x_e} 
		[\F_{\eph}^{\tmop{LW}} - \pf_{\Nmh}^e]
	\end{split} \label{eq:low.order.update}
\end{equation}
where $\pf_{\pph}^e = \pf ( \uep^n, \ueppone^n )$ is an
admissibility preserving finite volume numerical flux. Then, note that
\begin{equation}\label{eq:cell.avg.decomp}
	\avg{\uu}_e^{n + 1} = \sum_{p = 0}^N w_p  \atu_{e, p}^{n + 1} 
\end{equation}
Thus, if we can ensure that $\atu_{e, p}^{n + 1} \in \Uad$ for all
$p$, the scheme will be admissibility preserving in
means~\eqref{defn:mean.pres}. We do have $\atu_{e, p}^{n + 1} \in \Uad$ for $1 \le p \le N-1$ under appropriate CFL conditions because the finite volume fluxes are admissibility preserving. In
order to ensure that the updates $\atu_{e, 0}^{n + 1}, \atu_{e, N}^{n + 1}$ are also admissible, the flux limiting procedure of~{\cite{babbar2023admissibility}} is followed so that the high order numerical fluxes $\F_{e\pm\half}^{\tmop{LW}}$ are replaced by \textit{blended numerical fluxes} $\F_{e \pm \half}$. The procedure is explained here for completeness. We define an admissibility preserving lower order flux at the interface $\eph$
\[
\pf_{\eph} = \pf ( \uepoz^n, \ueN^n )
\]
Note that, for an RKFR scheme using Gauss-Legendre-Lobatto (GLL) solution
points, the definition of $\atu_{e, N}^{n + 1}$ will use
$\pf_{\eph}$ in place of $\F_{\eph}^{\tmop{LW}}$ and thus admissibility
preserving in means property will always be present. That is the argument
of~{\cite{Zhang2010b}} and here we demonstrate that the same argument can be applied to LWFR schemes by limiting $\F_{\eph}^\text{LW}$. We will explain the procedure for limiting $\F_{\eph}^{\tmop{LW}}$ to obtain $\F_{\eph}$; it
will be similar in the case of $\F_{\emh}$. Note that we want $\F_{\eph}$ to be such that the following are admissible
\begin{equation}
	\begin{split}
		\atu_0^{n + 1} & = \uepoz^n - \frac{\Delta t}{w_0 \Delta x_{e + 1}} 
		(\pf^{e + 1}_{\half} - \F_{\eph})\\
		\atu_N^{n + 1} & = \ueN^n - \frac{\Delta t}{w_N \Delta x_e}  (\F_{\eph}
		- \pf^e_{\Nmh})
	\end{split} \label{eq:low.order.tilde.update}
\end{equation}
We will exploit the admissibility preserving property of the finite volume fluxes to get
\begin{align*}
\utilow_0 & = \uepoz^n - \frac{\Delta t}{w_0 \Delta x_{e + 1}} 
(\pf^{e+1}_{\half} - \pf_{\eph}) \in \Uad\\
\utilow_N & = \ueN^n - \frac{\Delta t}{w_N \Delta x_e}  (\pf_{\eph} - \pf^e_{\Nmh}) \in \Uad
\end{align*}
Let $\{ \ad_k, 1 \le 1 \le K\}$ be the admissibility constraints~\eqref{eq:uad.form} of~\eqref{eq:con.law}. The time averaged flux is limited by iterating over the constraints. For each constraint, we can solve an optimization problem to find the largest
$\theta \in [0,1]$ satisfying
\begin{equation}
\ad_k ( \theta \atu_p^{n + 1} + (1 - \theta)  \utilow_p ) > \epsilon_p, \qquad p = 0, N \label{eq:optimization.problem}
\end{equation}
where $\epsilon_p$ is a tolerance, taken to be $\frac{1}{10}  \ad_k (\utilow_p)$~\cite{ramirez2021}. The optimization
problem is usually a polynomial equation in $\theta$. If $\ad_k$ is a concave
function of the conserved variables, we can
follow~{\cite{babbar2023admissibility}} and use the simpler but possibly sub-optimal approach of defining
\begin{equation}
	\theta = \min \left( \min_{p = 0, N} \frac{|\epsilon_p - \ad_k (\atu_p^{\text{low}, n + 1})|
		}{|\ad_k (\atu_p^{n + 1})  - \ad_k (\atu_p^{\text{low}, n + 1})| + \tmverbatim{eps}} , 1 \right) \label{eq:concave.theta}
\end{equation}
where $\tmverbatim{eps} = 10^{-13}$ is used to avoid a division by zero. In either case, by iterating over the admissibility constraints $\{ P_k \}$ of the
conservation law, the flux $\F_{\eph}^{\tmop{LW}}$ can be corrected as in Algorithm~\ref{alg:flux.correction}.
\begin{algorithm}
\caption{Flux limiting} \label{alg:flux.correction}
\begin{algorithmic}
\State $\F_{\eph} \leftarrow \F_{\eph}^{\tmop{LW}}$
\For{$k$=1:$K$}
\State $\epsilon_0, \epsilon_N \leftarrow \frac{1}{10}  \ad_k (\utilow_0), \frac{1}{10}  \ad_k (\utilow_N)$
\State $\tmop{Find} \theta$ by solving the optimization problem~\eqref{eq:optimization.problem} or by using~\eqref{eq:concave.theta} if $P_k$ is concave
\State $\F_{\eph} \leftarrow \theta \F_{\eph} + (1 - \theta) \pf_{\eph}$
\State $\atu_0^{n + 1} \leftarrow \uepoz^n - \frac{\Delta t}{w_0 \Delta x_{e + 1}}  (\pf^{e + 1}_{\half} - \F_{\eph})$
\State $\atu_N^{n + 1} \leftarrow \ueN^n - \frac{\Delta t}{w_N \Delta x_e}  (\F_{\eph} - \pf^e_{\Nmh})$
\EndFor
\end{algorithmic}
\end{algorithm}
After $K$ iterations, we will have $\ad_k (\atu_p^{n + 1}) \geq \epsilon_p$
for $p = 0, N$ and $1 \leq k \leq K$. Once this procedure is done at all the
interfaces, i.e., after Algorithm~\ref{alg:flux.correction} is performed, using the numerical flux $\F_{\eph}$ in the LWFR scheme will ensure that $\atu_{e, p}^{n +1}$~\eqref{eq:low.order.update} belongs to $\Uad$ for all $p$, implying $\au_e^{n + 1} \in \Uad$ by~\eqref{eq:cell.avg.decomp}. Then, we can use the scaling limiter
of~{\cite{Zhang2010b}} at all solution points and
obtain an admissibility preserving LWFR scheme for conservation laws in the absence of source terms. 
\subsection{Limiting time average sources} \label{sec:source.limiter}
After the flux limiting performed in Section~\ref{sec:flux.correction}, we will have an admissibility
preserving in means scheme~\eqref{defn:mean.pres} if the source term average $\avg{\bS}_e$ in~\eqref{eq:fravgup} is zero. In
order to get an admissibility preserving scheme in the presence of source terms, we will make a splitting of
the cell average update~\eqref{eq:fravgup}, which is similar to that of~\cite{meena2017}
\begin{equation}
	\au^{n + 1}_e = \frac{1}{2}  \left( \au^n_e - \frac{2 \Delta t}{\Delta x_e}  (\F_{\eph} - \F_{\emh})
	\right) + \frac{1}{2}  (\au^n_e + 2 \Delta t \avg{\bS^\text{LW}_e}
	) =: \au^{\F}_e + \au^{\bS^\text{LW}}_e \label{eq:F.S.split}
\end{equation}
where $\bS^{\tmop{LW}}_e$ denotes the time average source term in element $e$ computed with the approximate Lax-Wendroff procedure in Section~\ref{sec:approximate.lw}. With the
flux limiting performed in Section~\ref{sec:flux.correction}, we can ensure that cell average $\overline{\uu}^{\F}_e \in \Uad$ if twice the standard CFL is assumed (although, in the experiments we conducted, the CFL restriction used in~\cite{babbar2023admissibility} preserved admissibility). In order to enforce $\au^{\bS}_e \in \Uad$, $\bS^\text{LW}_e$ will be limited as follows. We will use the admissibility of the first order update using the source term
\[ 
\au^{\tmop{low}, n + 1}_e \assign \au^n_e + 2 \Delta t \avg{\bss}_e\in \Uad ,\qquad \avg{\bss}_e = \sum_{p=0}^N w_p \bss(\uep, \bx_{e,p}, t^n)
\]
which will be true under some problem dependent time step restrictions (e.g., Theorem 3.3.1 of~{\cite{Meena_Kumar_Chandrashekar_2017}}). Then, we will find a $\theta \in [0,1]$ so that for $\bS = \theta \bss+  (1-\theta) \bS^\text{LW}$, we will have $\avg{\uu}_e^{\bS} \in \Uad$. The $\theta$ can be found by iterating over admissibility constraints and using \eqref{eq:optimization.problem} or~\eqref{eq:concave.theta} where, for the first iteration, $\atu_p^{\text{low}, n + 1}$ is replaced by $\au^{\tmop{low}, n + 1}_e$ and $\atu_p^{n + 1}$ by $\au^{\bS^\text{LW}}_e$. Thus, a procedure analogous to Algorithm~\ref{alg:flux.correction} is used for limiting source terms. Then, replacing $\bS^\text{LW}$ by $\bS$
in~\eqref{eq:F.S.split}, we will have $\au^{\bS} \in \Uad$ and since $\F$ has
been corrected to ensure $\au^{\F}_e \in \Uad$, we will also have $\au^{n + 1}_e \in \Uad$. Thus, we have an admissibility preserving in means LWFR scheme~\eqref{defn:mean.pres} even in the presence of source terms. 
\section{Numerical results}\label{sec:numerical.results}
The numerical verification of admissibility preserving flux limiter~(Section~\ref{sec:flux.correction}) and admissibility of LWFR with source terms~(Section~\ref{sec:source.limiter}) is made through the Ten Moment equations~{\cite{Levermore_1996}} which we describe here. Here, the energy tensor is defined by the ideal equation of state $\ccE = \frac{1}{2}  \ccP + \frac{1}{2} \rho \bv \otimes \bv$ where $\rho$ is the density, $\bv$ is the velocity vector, $\ccR$ is the symmetric pressure tensor. Thus, we can define the 2-D conservation law with source terms
\[ 
\partial_t  \uu + \partial_{x_1}  \pf_1 + \partial_{x_2}  \pf_2 =
\bss^{x_1} \left( \uu \right) + \bss^{x_2} \left( \uu \right) 
\]
where $\uu =(\rho, \rho v_1, \rho v_2,  \cE_{11},  \cE_{12},  \cE_{22})$ and
\begin{equation}
	\pf_1 = \left[\begin{array}{c}
		\rho v_1\\
		\cR_{11} + \rho v_1^2\\
		\cR_{12} + \rho v_1 v_2\\
		\left( \cE_{11} + \cR_{11} \right) v_1\\
		\Eonetwo v_1 + \frac{1}{2}  ( \cR_{11} v_2 + \cR_{12} v_1 )\\
		\Etwotwo v_1 + \cR_{12} v_2
	\end{array}\right], \quad \pf_2 = \left[\begin{array}{c}
		\rho v_2\\
		\cR_{12} + \rho v_1 v_2\\
		\cR_{22} + \rho v_2^2\\
		\cE_{11} v_2 + \cR_{12} v_1\\
		\cE_{12} v_2 + \frac{1}{2}  ( \cR_{12} v_2 + \cR_{22} v_1 )\\
		\left( \cE_{22} + \cR_{22} \right) v_2
	\end{array}\right] \label{eq:tmp}
\end{equation}
The source terms are given by
\begin{equation}\label{eq:tenmom.source}
	\bss^{x_1} = \left[\begin{array}{c}
		0\\
		- \frac{1}{2} \rho \partial_x W\\
		0\\
		- \frac{1}{2} \rho v_1 \partial_x W\\
		- \frac{1}{4} \rho v_2 \partial_x W\\
		0
	\end{array}\right], \qquad  \bss^{x_2} = \left[\begin{array}{c}
		0\\
		0\\
		- \frac{1}{2} \rho \partial_y W\\
		0\\
		- \frac{1}{4} \rho v_1 \partial_y W\\
		- \frac{1}{2} \rho v_2 \partial_y W
	\end{array}\right] 
\end{equation}
where $W (x, y, t)$ is a given function, which models electron quiver energy in the laser~\cite{Berthon2015}. The admissibility set is given by
\[ 
\Uad = \left\{ \uu \in \re^6 | \rho \left( \uu \right) > 0, \quad \bx^T 
\ccP \left( \uu \right)  \bx > 0, \quad \bx \in \mathbb{R}^2 \backslash
\left\{ \bzero \right\} \right\} 
\]
which contains the states $\uu$ with positive density and positive definite
pressure tensor. The positive defininess of $\ccP$ implies that $\Poneone +
\Ptwotwo > 0$ and $\det \ccP = \Poneone  \Ptwotwo - \Ponetwo^2 > 0$. The
hyperbolicity of the system without source terms, along with its eigenvalues
are presented in Lemma 2.0.2 of~{\cite{Meena_Kumar_Chandrashekar_2017}}. The conditions for admissibiltiy preservation of the forward Euler method for the source terms, which are the basis for the source term limiting described in Section~\ref{sec:source.limiter}, are in Theorem 3.3.1 of~{\cite{Meena_Kumar_Chandrashekar_2017}}. All distinct numerical experiments from~\cite{Meena_Kumar_Chandrashekar_2017,meena2017,Meena2020} were performed and observed to validate the accuracy and robustness of the proposed scheme, but only some are shown here. The experiments were performed both with the TVB limiter used in~\cite{babbar2022} and the subcell-based blending scheme developed in~\cite{babbar2023admissibility}. As demonstrated in~\cite{babbar2023admissibility}, the subcell based limiter preserves small scale structures well compared to the TVB limiter. The use of TVB limiter is only made in this work to numerically validate that the flux limiting procedure of~Section~\ref{sec:flux.correction} preserves admissibility even without the subcells. The results shown are produced with TVB limiter unless specified otherwise.

The developments made in this work have been contributed to the package {\tt Tenkai.jl}~\cite{tenkai} written in {\tt Julia} and the setup files used for generating the results in this work are available in~\cite{icosahom2023_tmp}.
\subsection{Convergence test}
\begin{figure}
\centering
\begin{tabular}{cc}
\includegraphics[width=0.45\textwidth]{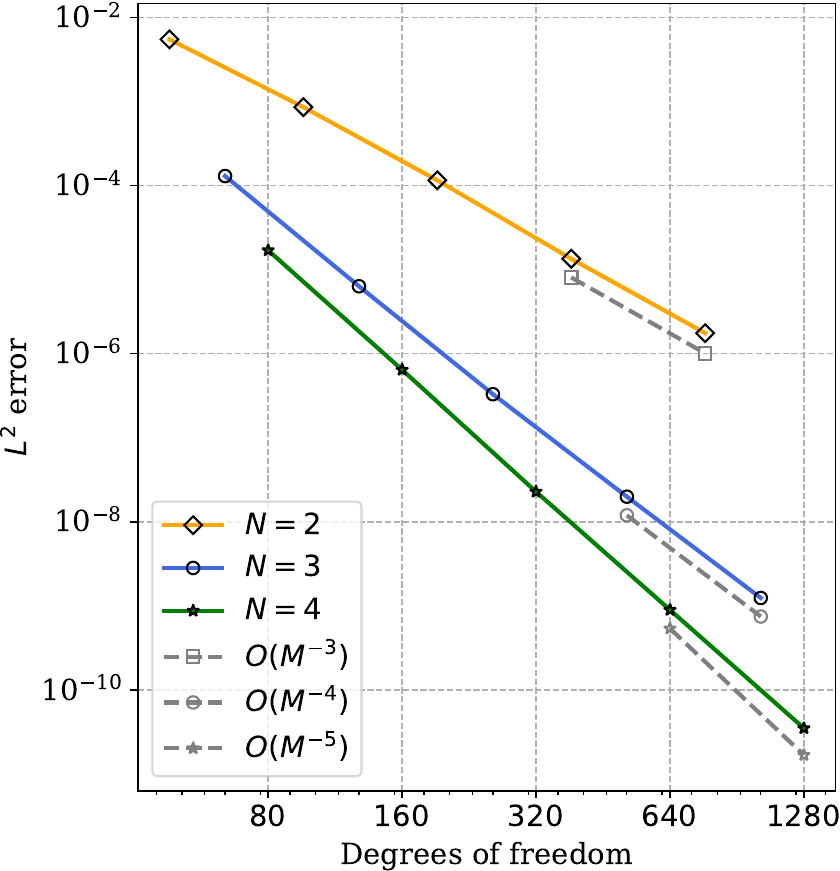} & \includegraphics[width=0.45\textwidth]{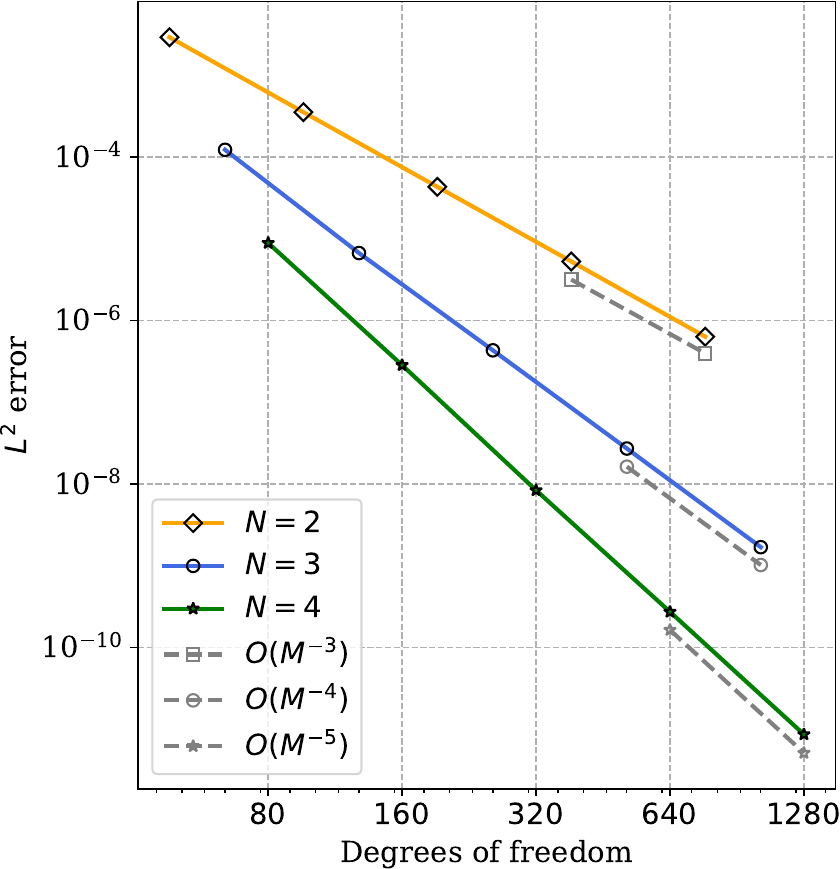} \\
(a) & (b)
\end{tabular}
\caption{Error convergence analysis of a smooth test with source terms for (a) $\rho$, (b) $\Poneone$ variable \label{fig:convergence}}
\end{figure}

This is a smooth convergence test from~{\cite{Biswas2021}} and requires no limiter. The domain is taken to be $\Omega = [-0.5, 0.5]$ and the potential for source terms~\eqref{eq:tenmom.source} is $W = \sin (2 \pi (x - t))$. With periodic boundary conditions, the exact solution is given by
\begin{equation*}
\begin{gathered}
{\rho}(x,t) =2+\sin  (2 {\pi}  (x-t)), \quad{\vone}(x,t) =1,\quad{\vtwo}(x,t) = 0\\
{\Poneone}=1.5+{\frac{1}{8}}  [\cos (4 {\pi}  (x-t)) -8 \sin (2 {\pi} 
(x-t))], \quad{\Ponetwo}(x,t) =0,\quad{\Ptwotwo}(x,t) =1
\end{gathered}
\end{equation*}
The solutions are computed at $t = 0.5$ and the convergence results for variable $\rho$ and $\Poneone$ are shown in Figure~\ref{fig:convergence} where optimal convergence rates are seen.
\subsection{Riemann problems}
\begin{figure}
\centering
\begin{tabular}{cc}
\includegraphics[width=0.45\linewidth]{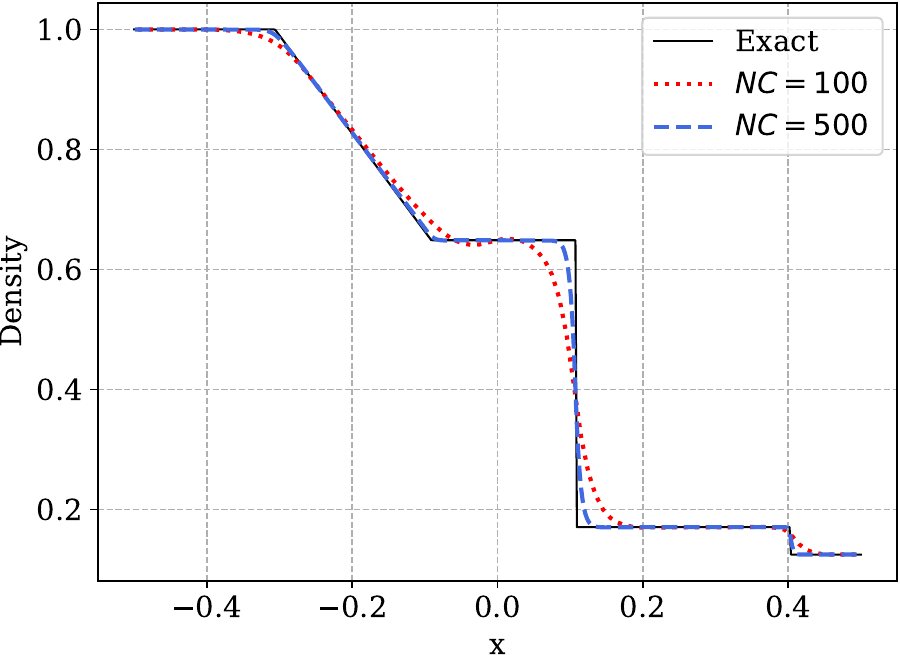} & \includegraphics[width=0.45\linewidth]{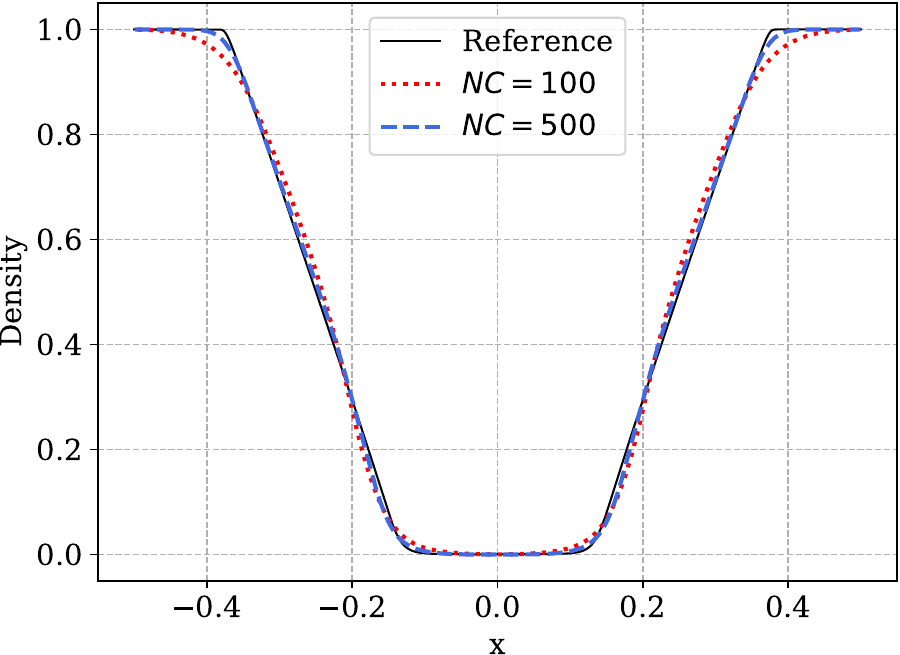} \\
(a)  & (b)
\end{tabular}
\caption{Density plots of numerical solutions with polynomial degree $N=2$ for (a) Sod's problem, (b) Two rarefaction (near vacuum) problem \label{fig:rp}}
\end{figure}

Here, we test the scheme on Riemann problems in the absence of source terms. The domain is $\Omega = [- 0.5, 0.5]$. The first problem is Sod's test
\[ 
( \rho, \vone, \vtwo, \Poneone, \Ponetwo, \Ptwotwo ) =
\begin{cases}
	(1, 0, 0, 2, 0.05, 0.6), \qquad & x < 0\\
	(0.125, 0, 0, 0.2, 0.1, 0.2),  & x > 0
\end{cases}
\]
The second is a problem from~\cite{Meena_Kumar_Chandrashekar_2017} with two rarefaction waves containing both low-density and low-pressure, leading to a near vacuum solution
\[ 
\left( \rho, v_1, v_2, \Poneone, \Ponetwo, \Ptwotwo \right) =
\begin{cases}
	(1, - 5, 0, 2, 0, 2), \qquad & x < 0\\
	(1, 5, 0, 2, 0, 2), & x > 0
\end{cases}
\]
The scheme is able to maintain admissibility in the near vacuum test and the results for both Riemann problems are shown in Figure~\ref{fig:rp} where convergence is seen under grid refinement.
\subsection{Shu-Osher test}
\begin{figure}
\begin{tabular}{cc}
\includegraphics[width=0.45\linewidth]{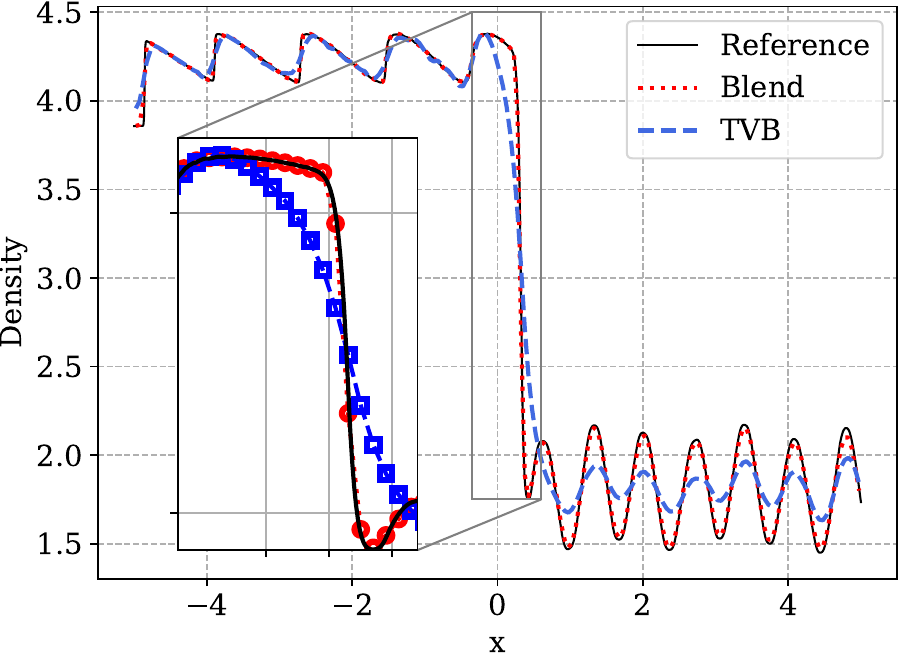} & \includegraphics[width=0.45\linewidth]{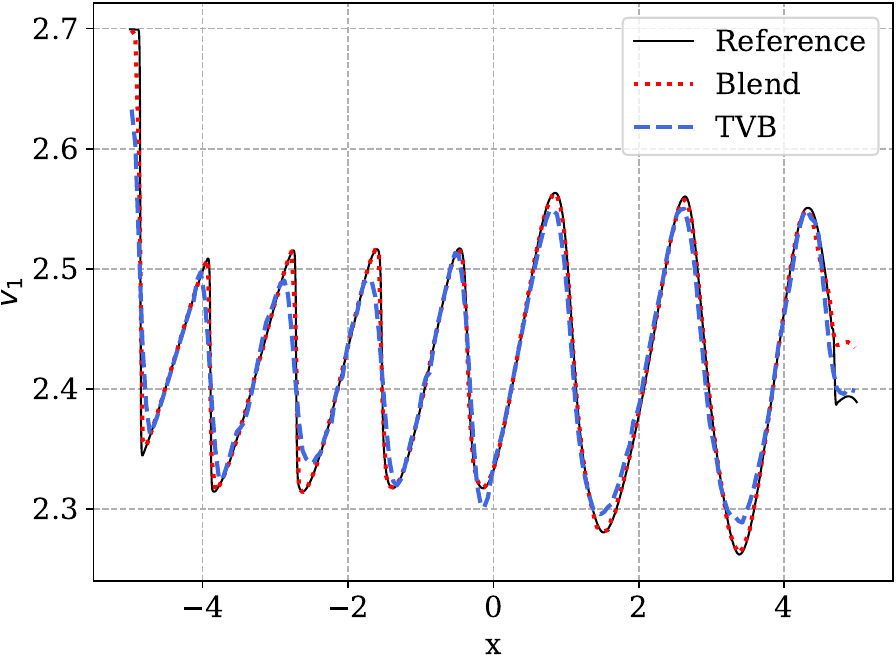} \\
(a) & (b)
\end{tabular}
\caption{Numerical solution for Shu-Osher problem with polynomial degree $N=4$ using TVB and blending limiter and we show (a) Density, (b) $v_1$ profiles. The density plot has an inset plot near the shock which compares number of cells smeared across the shock by blending and TVB limiter. \label{fig:shuosher}}
\end{figure}

\begin{figure}
\begin{tabular}{cc}
\includegraphics[width=0.45\linewidth]{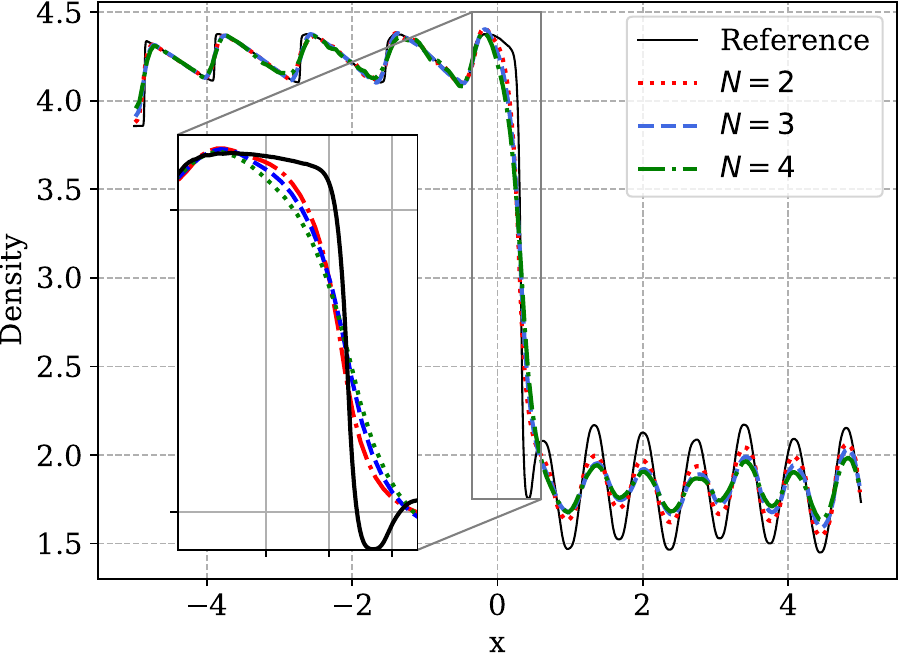} & \includegraphics[width=0.45\linewidth]{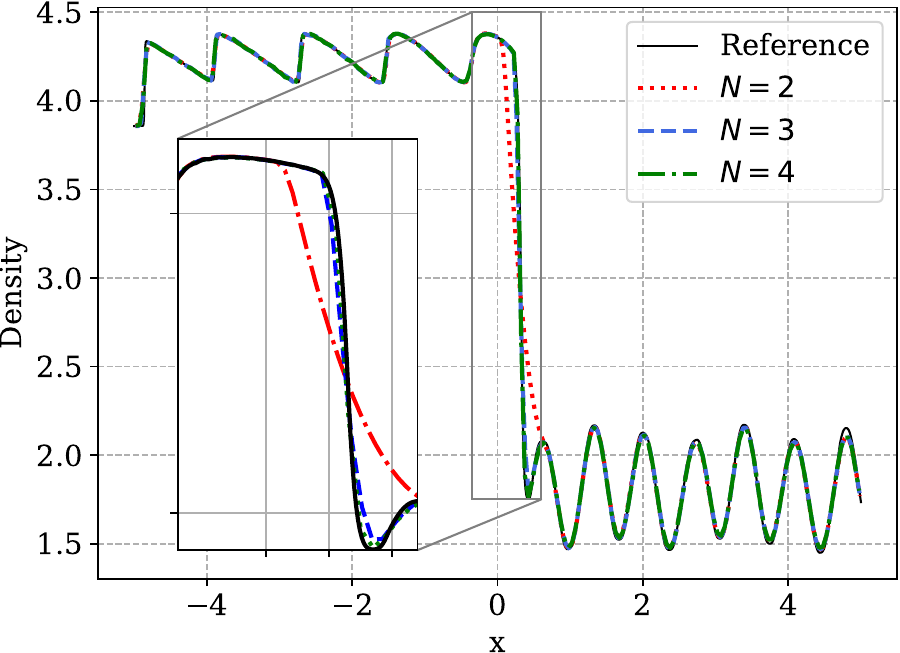} \\
(a) & (b)
\end{tabular}
\caption{Density profile of the numerical solution for Shu-Osher problem with using 1000 degrees of freedom (solution points) for polynomial degrees $N=2,3,4$ using (a) TVB and (b) blending limiter. \label{fig:shuosher.ndofs}}
\end{figure}

This is a modified version of the standard Shu-Osher test, taken from~\cite{Meena2020}. The primitive variables $\boldsymbol V = (\rho, v_1, v_2, \Poneone, \Ponetwo, \Ptwotwo)$ are initialized in domain $[-5,5]$ as
\[
\boldsymbol V = \begin{cases}
(3.857143, 2.699369, 0, 10.33333, 0, 10.33333),\qquad & \text{if } x \le -4 \\
(1 + 0.2 \sin(5x), 0, 0, 1, 0, 1), & \text{if } x > -4
\end{cases}
\]
The simulation is performed with polynomial degree $N=4$ using 200 elements and run till time $t = 1.8$ and the results with both blending and TVB limiter are shown in Figure~\ref{fig:shuosher} where, as expected, the blending limiter is giving much better resolution of the shock and high-frequency wave. In Figure~\ref{fig:shuosher.ndofs}, we show the numerical solutions for degrees $N=2,3,4$ using a grid with 1000 solution points (degrees of freedom) for each degree. Thus, the degree $N=4$ uses $200$ elements while degree $N=3$ uses $250$ elements. The numerical solutions in Figure~\ref{fig:shuosher.ndofs}a have been generated using the TVB limiter and it is seen that the accuracy actually degrades as we increase the degree. That is not the case for the results generated by blending limiter, shown in Figure~\ref{fig:shuosher.ndofs}b, where the accuracy improves with higher degrees. A great improvement is seen going from degree $N=2$ to $N=3$; $N=4$ is marginally better than $N=3$. This shows that the blending scheme, by preserving subcell information, is able to get the accuracy benefit of high order methods even in presence of shocks.
\subsection{Two dimensional near vacuum test}
\begin{figure}
\centering
\begin{tabular}{cc}
\includegraphics[width=0.35\linewidth]{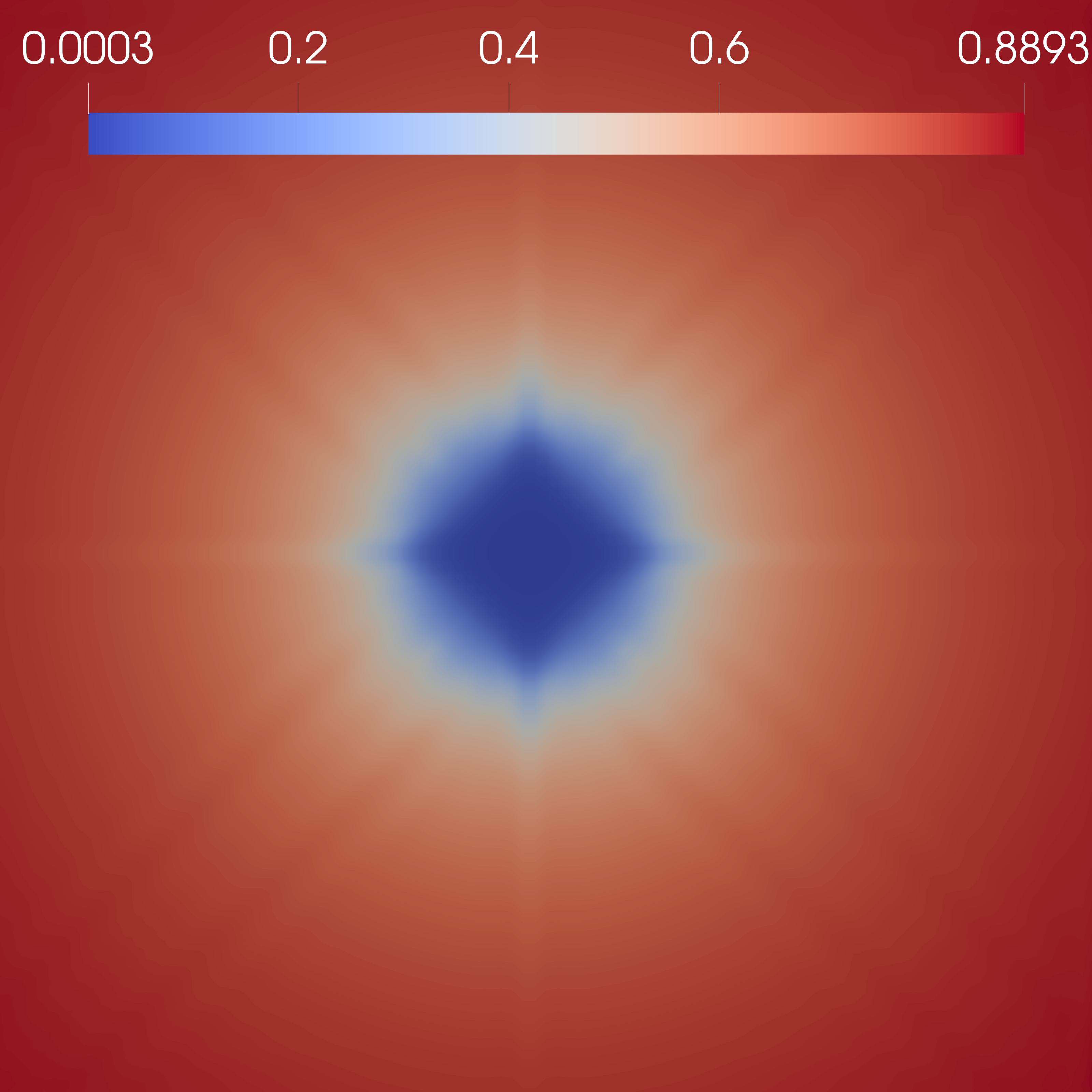} & \includegraphics[width=0.48\linewidth]{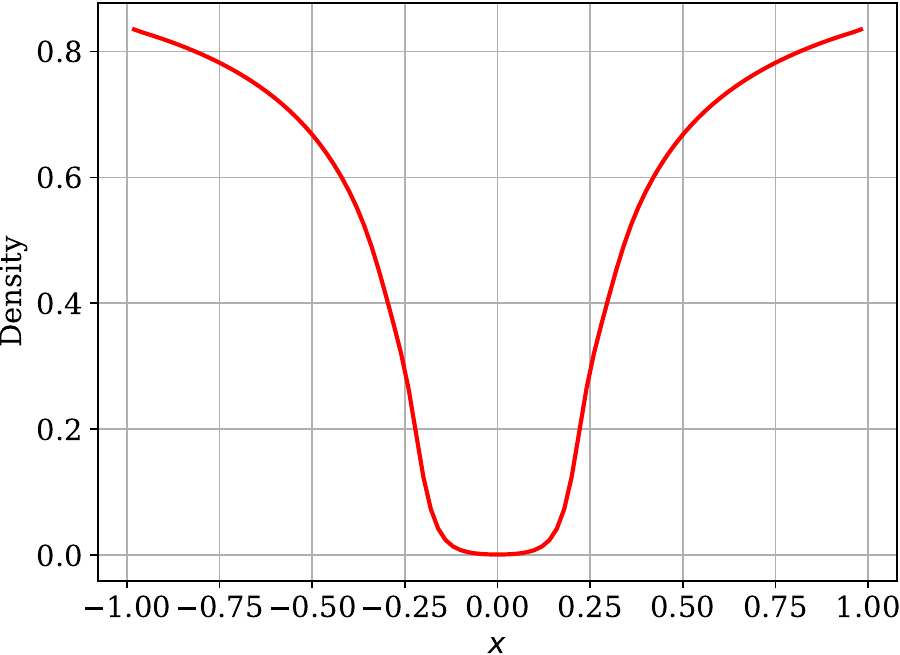} \\
(a) & (b)
\end{tabular}
\caption{2-D near vacuum test. Density plot of numerical solution with degree $N=2$ on a $100^2$ element mesh (a) Pseudocolor plot (b) Solution cut along the line $y=0$}
\label{fig:2d.near.vacuum}
\end{figure}

This is a near vacuum test taken from~\cite{Meena_Kumar_Chandrashekar_2017}, and is thus another verification of our admissibility preserving framework. The domain is $\Omega = [- 1, 1]^2$ with outflow boundary conditions. The initial conditions are
\begin{equation*}
\begin{gathered}
{\rho}=1,\quad{\Poneone}={\Ptwotwo}=1,\quad{\Ponetwo}=0, \quad
{\vone}=8 f\tmrsub{s}(r)  \cos {\theta}, \quad{\vtwo}=8 f\tmrsub{s}(r)  \sin
{\theta}
\end{gathered}
\end{equation*}
where $r = \sqrt{x^2 + y^2}$, $\theta = \arctan (y / x) \in [- \pi, \pi]$ and
$s = 0.06 \mathLaplace x$ for mesh size $\mathLaplace x(=\Delta y)$ of
the uniform mesh. The $f_s (r)$ smoothens the velocity profile near the origin
as $\theta$ is not defined there
\[ 
f_s (r) = \begin{cases}
	- 2 \left( \frac{r}{s} \right)^3 + 3 \left( \frac{r}{s} \right)^2, \quad
	& \tmop{if} r < s\\
	1, & \text{otherwise}
\end{cases}
\]
The numerical solution computed using polynomial degree $N=2$ and 100 elements is shown at the time $t = 0.02$. The results are shown in Figure~\ref{fig:2d.near.vacuum} and are similar to those seen in the literature.
\subsection{Realistic simulation}
\begin{figure}
\centering
\begin{tabular}{cc}
\includegraphics[width=0.345\linewidth]{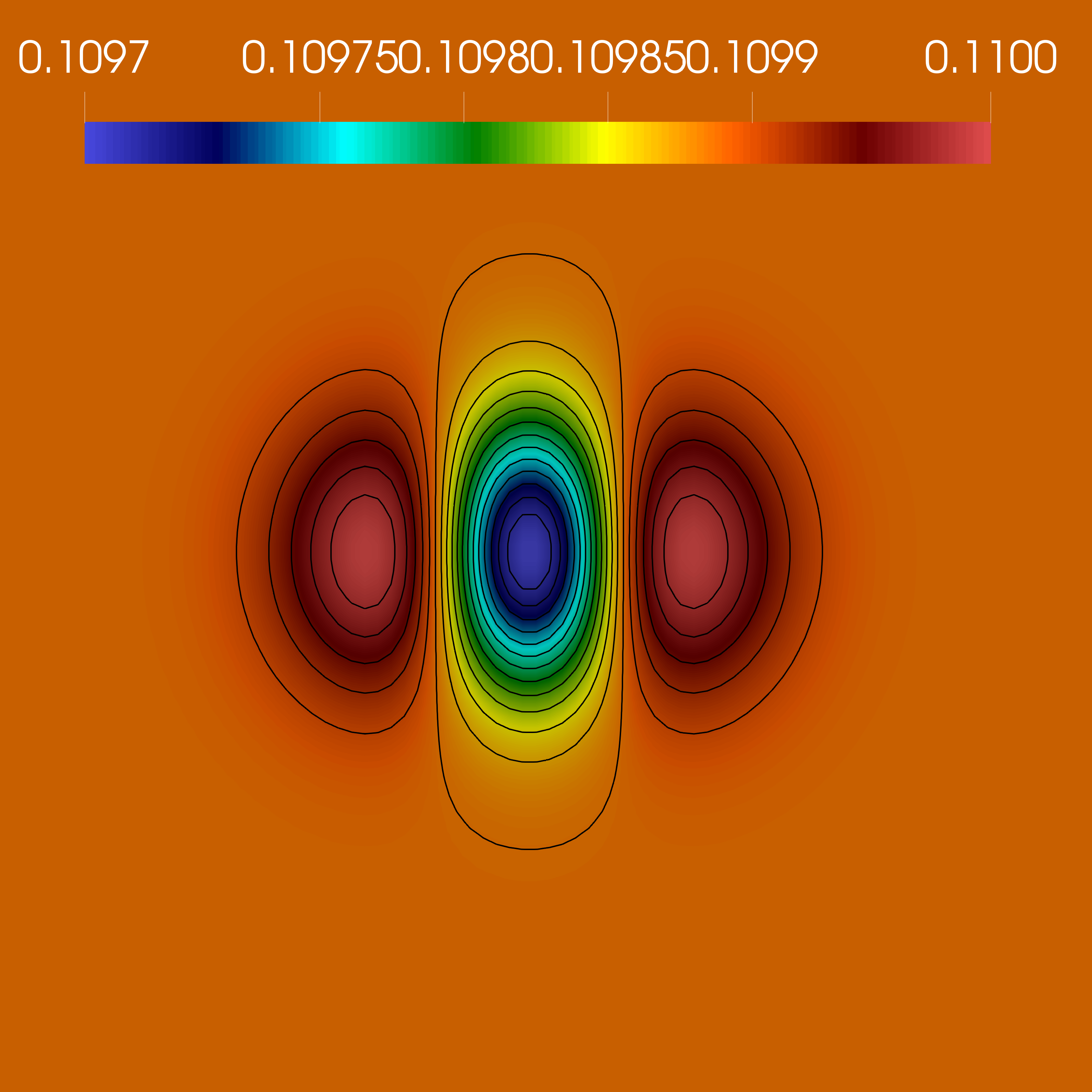} & \includegraphics[width=0.46\linewidth]{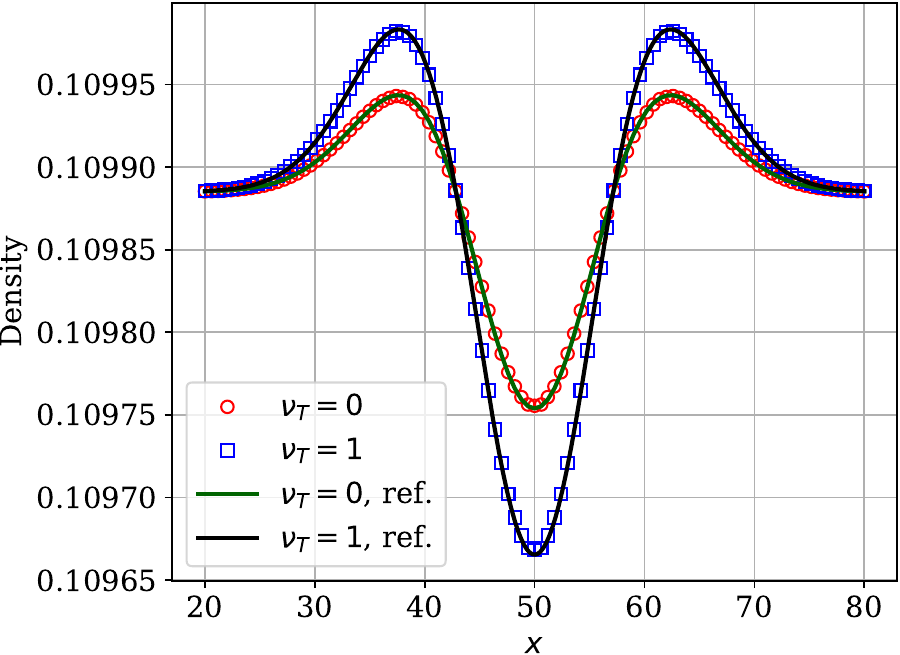} \\
(a) & (b)
\end{tabular}
\caption{Realistic simulation. Density profile computed with degree $N=4$ on $100^2$ element mesh. (a) Pseudocolor color plot (b) Cut at $y=4$ comparing different absorption coefficient $\nu_T$. \label{fig:realistic}}
\end{figure}

Consider the domain $\Omega = [0, 100]^2$ with outflow boundary conditions. The uniform initial condition is taken to be
\[ 
\rho = 0.109885, \quad \vone = \vtwo = 0, \quad \Poneone = \Ptwotwo = 1,
\quad \Ponetwo = 0 \]
with the electron quiver energy $W (x, y, t) = \exp ( - 0.01 (( x - 50 )^2 + (y - 50 )^2) )$. The source term is taken from~{\cite{Berthon2015}}, and only has the $x$ component, i.e, $\bss^y \left( \uu \right) = \bzero$, even though $W$ continues to depend on $x$ and $y$. An additional source corresponding to energy components $\bss_E = (0, 0, 0, \nu_T \rho W, 0, 0)$ is also added where $\nu_T$ is an absorption coefficient. Thus, the source terms are $\bss = \bss_x + \bss_E$. The simulation is run till $t = 0.5$ on a grid of 100 cells. The blending limiter from~\cite{babbar2023admissibility} was used in this test as it captured the smooth extrema better. The density plot with a cut at $y = 4$ is shown in Figure~\ref{fig:realistic}.
\section{Conclusions} \label{sec:conclusion}
A generalized framework was developed for high order admissibility preserving
Lax-Wendroff (LW) schemes. The framework is a generalization of~\cite{babbar2023admissibility} as it is independent of the shock capturing scheme used and can, in particular, be used without the subcell based limiter of~\cite{babbar2023admissibility}. The framework was also shown to be an extension
of~{\cite{Zhang2010b}} to LW. The LW scheme was extended to be applicable to
problems with source terms while maintaining high order accuracy. Provable
admissibility preservation in presence of source terms was also obtained by
limiting the time average sources. The claims were numerically verified on the Ten Moment
problem model where the scheme showed high order accuracy and robustness.
\bibliographystyle{spbasic}
\bibliographystyle{siam}
\bibliography{references.bib}
\end{document}